\documentclass[a4paper]{amsart}
\usepackage{graphicx}
\usepackage{hyperref}
\hypersetup{
    pdftitle={Theorem of Kerekjarto}, 
    pdfauthor={A. Constantin and B. Kolev} 
    }
\newtheorem{thm}{Theorem}[section]
\newtheorem{cor}[thm]{Corollary}
\newtheorem{lem}[thm]{Lemma}
\newtheorem{prop}[thm]{Proposition}
\theoremstyle{definition}

\theoremstyle{remark}
\newtheorem{rem}[thm]{Remark}
\numberwithin{equation}{section}

\begin{document}

\title[Theorem of K{\'e}r{\'e}kjart{\'o}]{The theorem of K{\'e}r{\'e}kjart{\'o} on periodic
homeomorphisms of the disc and the sphere}%

\author[A. Constantin]{Adrian CONSTANTIN}%
\address{Lund University, Dept. of Maths. P.O  Box 118, S-22100 Lund, Sweden}%
\email{adrian@maths.lth.se}%

\author[B. Kolev]{Boris Kolev}%
\address{CMI, 39, rue F. Joliot-Curie, 13453 Marseille cedex 13, France}%
\email{kolev@cmi.univ-mrs.fr}%

\thanks{The authors express their gratitude to J{\'e}r{\^o}me Fehrenbach, Lucien Guillou and
Toby Hall for the discussions that helped to improve this paper.}%
\subjclass{55M99 (54H20)}%

\date{1994}

\begin{abstract}
We give a modern exposition and an elementary proof of the
topological equivalence between periodic homeomorphisms of the
disc and the sphere and euclidean isometries.
\end{abstract}

\maketitle

\section{Introduction}

In 1919, K{\'e}r{\'e}kjart{\`o} published the first proof of the topological
equivalence between periodic homeomorphisms of the disc and the
sphere and euclidean isometries \cite{Kerekjarto1}. In the same
journal just following K{\'e}r{\'e}kjart{\`o}'s article, Brouwer
\cite{Brouwer} gave his own argument for these theorems,
explaining that these results had been known to him for a long
time and that they were consequences of some earlier and slightly
different theorems of his on periodic homeomorphisms of compact
surfaces. However, Brouwer's proof is not easy to follow and the
proof of K{\'e}r{\'e}kjart{\`o} was just sketched and contained a gap.

It was only in 1934 that a complete proof of this important
theorem was presented by Eilenberg \cite{Eilenberg}. More recently
Epstein \cite{Epstein} has reconsidered the question for pointwise
periodic homeomorphisms (each point is periodic under f but the
period $n(x)$ depends on $x$ and may not be bounded). Because of
the importance of these results and since no modern exposition of
them seems to be found in the literature, the authors have thought
that it would be useful to present a modern and elementary proof.
The essential arguments, however, remain those of
\cite{Brouwer,Eilenberg,Kerekjarto1}.


\section{Background and Definitions}

Let $X$ be a topological space and $f$ a homeomorphism of $X$. We
say that $f$ is periodic if there is an integer $n > 0$ such that
$f^{n} = Id$. The period of $f$ is the smallest positive integer
$n$ with this property.

As we will use them without further justifications, let us first
recall some basic properties of one- dimensional maps.

Let $f : I\to I$ be a periodic homeomorphism of the unit interval.
If $f$ preserves the endpoints then $f$ is the identity map. If
$f$ exchanges the endpoints then $f^{2}= Id$ and $f$ is conjugate
to the reflection map $x\mapsto  1 - x$. Similarly, a periodic
homeomorphism of the real line $\mathbb{R}$ is the identity map or
is a conjugate of the involution $x\mapsto 1-x$ according to
whether it is an increasing or a decreasing function.

Let $f:S^{1}\to S^{1}$ be a periodic homeomorphism of period $n$
of the unit circle. If $f$ is order-preserving then the rotation
number of $f$, $\rho (f) = k/n$, where $k$ and $n$ are coprime
(see \cite{Devaney} for an excellent exposition on rotation
numbers) and $f$ is conjugate to a rotation of angle $2k\pi/n$. If
$f$ is order-reversing then $f$ has exactly two fixed points,
$f^{2}$ is the identity map and the two arcs delimited on $S^{1}$
by the fixed points of $f$ are permuted by $f$.

A metric space $X$ is path connected if there exists a continuous
map from the unit interval $[0, 1]$ into $X$ which joins any two
given points. It is arcwise connected if there is a topological
embedding of $[0, 1]$ into $X$ which joins any two given distinct
points. In fact, it can be shown that the two notions are
equivalent (see \cite[Theorem 4.1]{Whyburn} or \cite[Lemma
16.3]{Milnor}).

\begin{lem}\label{2lem1}
A metric space X is path connected if and only if it is arcwise
connected.
\end{lem}

A useful characterization of path connected spaces is given in
term of local connectivity. A metric space $X$ is locally
connected if each point of $X$ possesses arbitrary small connected
neighborhoods. The following can be shown (see \cite[Theorem
3.15]{Hocking+Young} or \cite[Lemma 16.4]{Milnor}):

\begin{lem}\label{2lem2}
A compact, connected and locally connected metric space is
pathwise connected.
\end{lem}

Another important ingredient used in this article, and in fact the
ultimate result we will need, is the famous Jordan-Schoenflies
theorem on simple closed curves in the plane (see
\cite{Cairns,Maehara} or \cite[Theorem 17.1]{Newman}).

\begin{thm}[Jordan-Schoenflies]\label{Jordan}
Every simple closed curve J divides the plane into exactly two
components of each of which it is the complete boundary and the
closure of the bounded component can be mapped topologically onto
the closed unit disc.
\end{thm}

In what follows, a \textsl{closed} topological disc (or just a
topological disc) $D$ is the image under a topological embedding
of the closed unit disc and we write $D^{\circ}$ for its interior
and $\partial D$ for its boundary. However, the closure of a
bounded open set which is homeomorphic to the open unit disc is
not necessarily a closed topological disc \cite[Chapter
15]{Milnor}.

\begin{prop}\label{2prop4}
Let $D_{1}, D_{2},\ldots, D_{n}$ be a finite number of closed
topological discs in the plane and $J^{\circ}$ be any connected
component of $\cap_{i=1}^{n}D_{i}^{\circ}$. Then $\partial J$ is a
simple closed curve and J, the closure of $J^{\circ}$ is a
topological disc.
\end{prop}

\begin{proof}
We will use induction on $n$, the number
of discs. If $n = 1$ this is just the Jordan- Schoenflies theorem,
so let us suppose that the result holds for some $n(n\geq 1)$ and
let $J^{\circ}$ be any component of the complement of $n + 1$
topological discs $D_{1}, D_{2},\ldots, D_{n+1}$ in the plane. Let
$K^{\circ}$ be the component of $\cap^{n}_{i=1} D^{\circ}_{i}$
that contains $J^{\circ}$. By induction, its closure $K$ is a
topological disc. Since $J^{\circ}$ is a component of
$K^{\circ}\cap D^{\circ}_{n+1}$, it suffices to show that the
result holds for two discs $D_{1}$ and $D_{2}$ (see Figure
\ref{fig1}). Set $C_{i} = \partial D_{i}$ for $i = 1, 2$ and let
$J$ be the closure of a component of $D_{1}^{\circ}\cap
D_{2}^{\circ}$. We have that $\partial J \neq \emptyset$ and
$\partial J \subset C_{1}\cup C_{2}$. If $\partial J$ is entirely
contained in one of the two curves, say $C_{1}$, then $J = D_{1}$
and the lemma is proved. We can thus suppose that $\partial J
\not\subset C_{1}$ and $\partial J \not\subset C_{2}$.

Let $x\in \partial J$, $x \not\in C_{2}$. Then $x \in C_{1}\cap
D_{2}^{\circ}$, and we can find an arc $\gamma$ in $C_{1}$ such
that:
\begin{equation}
    x \in \gamma ,\quad  \gamma\subset\partial J
    ,\quad \gamma\setminus\partial\gamma \subset D_{2}^{\circ} ,\quad \partial\gamma
    \subset C_{2}.
\end{equation}

\begin{figure}
\begin{center}
\includegraphics[scale=0.5]{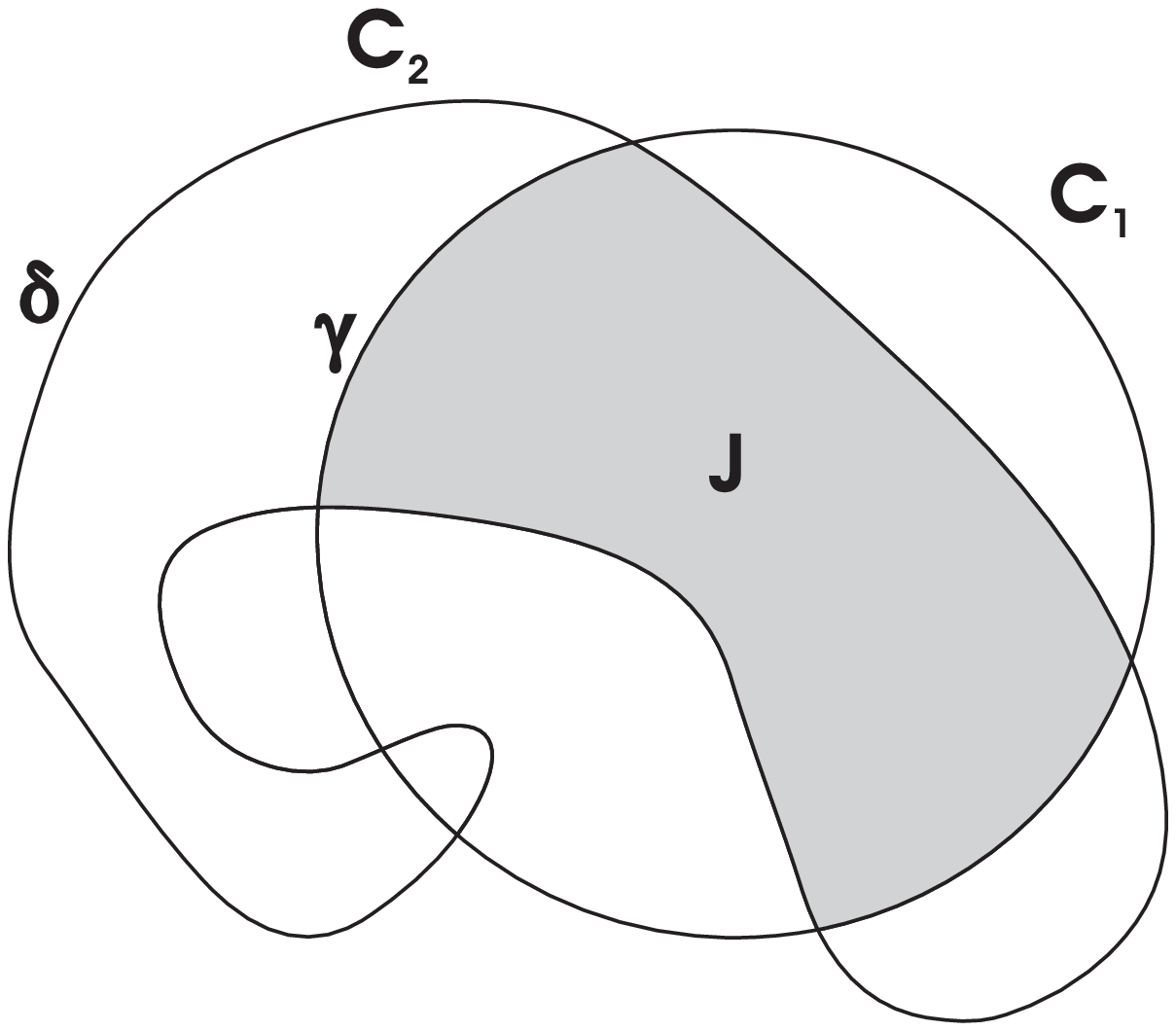}
\caption{} \label{fig1}
\end{center}
\end{figure}

The endpoints of $\gamma$ determine on $C_{2}$ an arc $\delta$
disjoint from $J^{\circ}$ and such that $\delta \cap J =
\partial\delta$. We note that there is an at most countable
family of such arcs $\gamma $, noted $ (\gamma_{i})_{i\in
\mathbb{N}}$ and that $diam(\gamma_{i})\to 0$ as $i\to\infty$. The
boundary of $J$ is the simple closed curve obtained from $C_{2}$
when substituting the arcs $\gamma_{i}$ for the arcs $\delta_{i}$
and $J$ is a topological disc by the Jordan-Schoenflies
theorem.
\end{proof}

The following remarkable property of periodic homeomorphisms which
is a direct consequence of \ref{2prop4} is true in a more general
setting than the plane $\mathbb{R}^{2}$, namely in topological
manifolds of dimension $2$ because of its local nature. We will
give it in that context since we will use it for the disc and the
sphere, repeatedly in this article.

\begin{lem}\label{2lem5}
Let $f: S \to S$   be a periodic homeomorphism of an arbitrary
$2$-dimensional topological manifold $S$ and let $x\in Fix(f)$, a
fixed point of $f$. Then for any neighbourhood $N$ of $x$, there
exists a topological disc $\Delta_{x}$ such that:
\begin{enumerate}
    \item $\Delta_{x}\subset N$,
    \item $\Delta_{x}$ is a neighbourhood of x,
    \item $f(\Delta_{x}) = \Delta_{x}$.
\end{enumerate}
\end{lem}

\begin{proof}
We can first assume that $N$ and its image
under $f$, $f(N)$, are contained in some local chart $U$
homeomorphic with $\mathbb{R}^{2}$ and we will continue to call
$x$ and $N$ the corresponding point and set in $\mathbb{R}^{2}$.
Let $D_{x}$ be an euclidean disc of centre $x$ and radius $\eta$
where $\eta> 0$ is chosen such that $f^{k}(D_{x})\subset N$ for $k
= 0, 1,\ldots, n - 1$ and let $C_{x}$ be its boundary. Let
$\Delta_{x}$ be the closure of the component of the invariant set
$\cap^{n-1}_{k=0}f^{k}(D^{\circ}_{x})$ which contains $x$. By
\ref{2prop4}, $\Delta_{x}$ is a topological disc which is
invariant under $f$ (components are sent to components by a
homeomorphism) and satisfies the three assertions of the lemma.
\end{proof}

\begin{rem}\label{2rem6}
The boundary $\gamma_{x}$ of $\Delta_{x}$ which
is an invariant simple closed curve, is contained in
$\cup_{k=0}^{n-1}f^{k}(C_{x})$.
\end{rem}


\section{Periodic Homeomorphisms of the Disc}

\begin{thm}\label{3thm1}
Let $f : D^{2} \to D^{2}$  be a periodic homeomorphism. Then there
exists $r \in O(2)$ and a homeomorphism $h : D^{2} \to D^{2}$ such
that $f=hrh^{-1}$.
\end{thm}

Before attacking the proof of the result above, let us first look
at a special case of Theorem \ref{3thm1}, namely:

\begin{prop}\label{3prop2}
Let $f : D^{2} \to D^{2}$ be a periodic homeomorphism such that
$f_{/\partial D^{2}} = Id$. Then $f = Id$.
\end{prop}

\begin{proof}
Let $d$ be an arbitrary diameter of $D^{2}$
with endpoints $A$ and $B$ and let $\Delta$ be one of the two
connected components of $D^{2} \setminus d$. The set:
\begin{equation}
    E = \bigcap_{i=1}^{n}f^{i}(\Delta^{\circ})
\end{equation}

is invariant under $f$ and the closure of each of its components
is a topological disc.

\begin{figure}
\begin{center}
\includegraphics[scale=0.5]{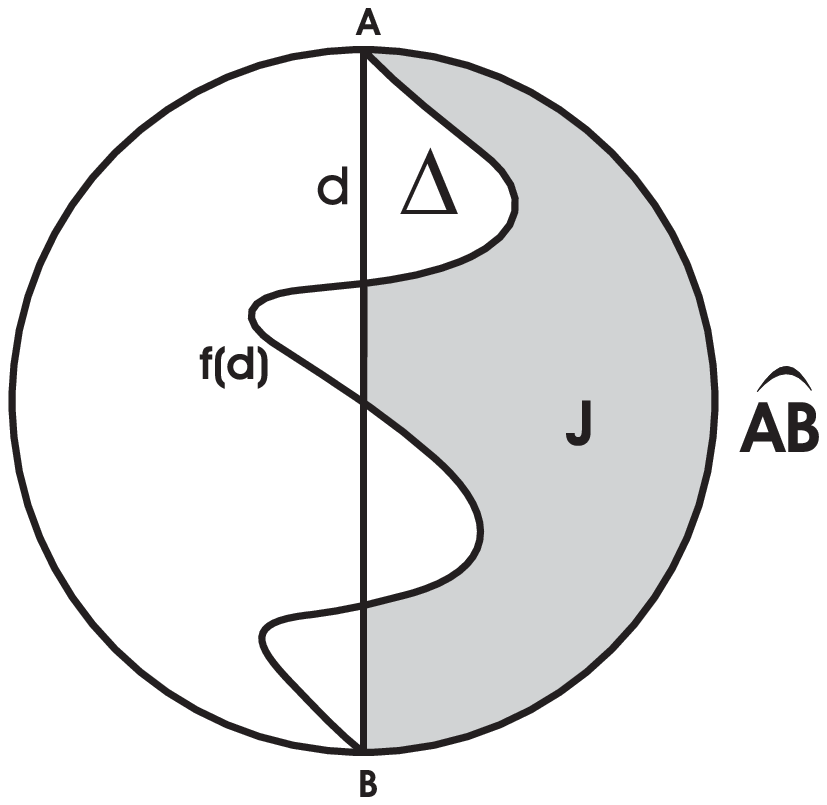}
\caption{} \label{fig2}
\end{center}
\end{figure}

Let $\widehat{AB}$ be the arc of circle joining $A$ to $B$ in the
boundary of $\Delta$. Since $f^{i}(\widehat{AB}) = \widehat{AB}$
for all $i$, there exists a component of $E$, say $J^{\circ}$,
whose closure $J$ contains $\widehat{AB}$(see Figure \ref{fig2}).
By \ref{2prop4}, $J$ is a topological disc which is invariant
under $f$. We can write $\partial J = \widehat{AB}\cup\delta$
where $\delta$ is an $f$-invariant, simple arc with endpoints $A$
and $B$ such that:
\begin{equation}
    \delta\subset\bigcup_{i=1}^{n}f^{i}(d).
\end{equation}

Since $f(A) = A$ and $f(B) = B$, $f_{/\delta}= Id$. Let $x$ be a
point of the arc $\delta$. There  exists  $i \in \{ 1,\ldots, n
\}$ such that $x \in f^{i}(d)$ and $x = f^{n-i}(x)\in d$ so that
$\delta = d$ and $f_{/d} = Id$. Since the diameter $d$ was chosen
arbitrarily, we have shown that $f = Id$ on $D^{2}$.
\end{proof}

From now on, $f$ will denote a periodic homeomorphism of the disc
of period $n$ with $n > 1$. In the sequel of this section, we
prove Theorem \ref{3thm1}, first investigating the structure of
the fixed point set of $f$.

\begin{prop}\label{3prop3}
Suppose $f : D^{2}\to D^{2}$ is a periodic homeomorphism of period
$n (n > 1)$, then:
\begin{enumerate}
    \item if $f$ is orientation-preserving, $Fix(f)$ is reduced to a single point which is not on
the boundary of $D^{2}$ and for $1 \leq i \leq n - 1$,
$Fix(f^{i})=Fix(f)$,
    \item if $f$ is orientation-reversing,  $f^{2} = Id$ and $Fix(f)$ is a simple arc
which divides $D^{2}$ into two topological discs which are
permuted by $f$.
\end{enumerate}
\end{prop}

\begin{proof}
Suppose first that $f$ is
orientation-preserving. By Brouwer fixed point theorem, $f$ has at
least one fixed point. Since $f_{/\partial D^{2}}$ is
orientation-preserving and periodic, $f$ has no fixed point on
$\partial D^{2}$. Otherwise $f$ would be the the identity map on
$\partial D^{2}$ and using \ref{3prop2}, $f$ would be the identity
map on the whole disc which is excluded by hypothesis. Therefore,
$f$ has at least one fixed point in $D^{2}\setminus
\partial D^{2}$ which we can assume to be, up to conjugacy, $O$,
the center of the disc.

Let $A = D^{2}\setminus{\{ O \}}$. $A$ is a half open annulus
which is invariant under $f$. Suppose now that an iterate $f^{i}$
of $f$ has a fixed point $x_{0}\in A$. Let $\widetilde{x_{0}}$ be
a lift of $x_{0}$ to the universal covering space $\widetilde{A}$
of $A$ and $G$ be the lift of $f^{i}$ such that
$G(\widetilde{x_{0}}) = x_{0}$. $G^{n}$ is a lift of $Id$ which
fixes one point, thus $G^{n} = Id$. In particular, $G_{/\partial
\widetilde{A}}$ is a periodic and orientation-preserving
homeomorphism of the line, thus $G = Id$ on
$\partial\widetilde{A}$. Therefore, $f^{i}= Id$ on $\partial
D^{2}$ and, according to \ref{3prop2}, $f^{i} = Id$ on the whole
disc, so that $i$ is a multiple of $n$ according to the definition
of $n$.

Suppose now that $f$ is orientation-reversing. In that case, $f$
has exactly two fixed points on $\partial D^{2}$ which we denote
by $A$ and $B$ and $f^{2}$ is the identity map on $\partial
D^{2}$, therefore, by \ref{3prop2}, $f^{2} = Id$ on $D^{2}$.

We assert that $Fix(f)$ is connected. For if not, we can find two
nonempty compact sets $K_{1}$, and $K_{2}$ such that:
\begin{equation}
   Fix(f) = K_{1}\cup K_{2}, \qquad K_{1}\cap K_{2}=\emptyset
\end{equation}

If $A \in K_{1}$, and $B \in K_{2}$, it is then possible to
construct a simple arc $\gamma$ in $D^{2} \setminus (K_{1}\cup
K_{2})$ which intersect $\partial D^{2}$ only on its endpoints and
which separates $A$ from $B$. Using the same argument as the one
used in the proof of \ref{3prop2}, we can show the existence of an
$f$-invariant simple arc:
\begin{equation}
    \delta\subset\bigcup_{i=0}^{n-1}f^{i}(\gamma )\subset
    D^{2}\setminus Fix (f)
\end{equation}

which separates $A$ from $B$. But $f$ must then have a fixed point
on $\delta$ which gives a contradiction. Therefore we can suppose
that one of the two compact sets, say $K_{1}$, is contained in
$D^{2}\setminus\partial D^{2}$. In that case, it is possible to
construct a simple closed curve $c\subset D^{2}\setminus\partial
D^{2}$ which does not meet $K_{1}\cup K_{2}$ and such that the
topological disc it bounds contains at least one point of $K_{1}$.
Using similar arguments as those of the proof of \ref{2lem5}, we
can find an $f$-invariant topological disc in
$D^{2}\setminus\partial D^{2}$ whose boundary contains no fixed
point. This gives again a contradiction, since any simple closed
curve which bounds an invariant disc has exactly two fixed points
of $f$.

The previous arguments applied to an arbitrarily small invariant
topological disc around a fixed point given by \ref{2lem5} shows
that $Fix(f)$ is also locally connected and by \ref{2lem2},
$Fix(f)$ is therefore pathwise connected. In view of \ref{2lem1},
there exists a simple arc $\gamma$ in $Fix(f)$ which joins $A$ and
$B$. This arc divides $D^{2}$ into two topological discs
$\Delta_{1}$ and $\Delta_{2}$ by the Jordan-Schoenflies theorem.
$D^{2}\setminus\gamma$ is obviously invariant under $f$ and the
two arcs on $\partial D^{2}$ delimited by $A$ and $B$ are permuted
by $f$, therefore $f(\Delta_{1}) = \Delta_{2}$, $f(\Delta_{2}) =
\Delta_{1}$ and $Fix(f)$ is reduced to $\gamma$.
\end{proof}

\begin{proof}[Proof of Theorem \ref{3thm1}]
Suppose first that $f$ is orientation-preserving. By \ref{3prop3},
we can suppose that $Fix(f) = \{ O \}$, the center of the disc.
Since $f_{/\partial D^{2}}$ is a periodic homeomorphism of period
$n$, the rotation number of $f_{/\partial D^{2}}$, $\rho
(f_{/\partial D^{2}}) = k / n$ where $k$ and $n$ are coprime. We
are going to prove that $f$ is conjugate to a rotation by angle
$2k\pi / n$ around the origin. Without loss of generality, we can
assume that $k = 1$. Indeed, suppose the result holds if $\rho
(f_{/\partial D^{2}}) = 1 / n$. Then, if $k > 1$ we replace $f$ by
$f^{j}$ where $j\in \mathbb{N} $ is such that $jk \equiv 1 (mod
n)$. Then  $\rho (f^{j}_{/\partial D^{2}}) = 1 / n$, thus $f^{j}$
is conjugate to a rotation by angle $2\pi /n$ around the origin
and since $(f^{j})^{k}=f$, it follows that $f$ is conjugate to a
rotation by angle $2k\pi /n$.

Let us consider the quotient space $D^{2}/f$ where two points are
identified if they belong to the same orbit under $f$. $D^{2}/f$
is endowed with the quotient topology. It is a compact and
pathwise connected metric space, the metric being defined by:
\begin{equation}
    d(\pi(x),\pi(y)) = \inf_{0\leq h,k\leq n-1}
    \{d(f^{k}(x),f^{k}(y))\},
\end{equation}
where $\pi : D^{2} \to D^{2}/f$ is the canonical projection.

By \ref{2lem1}, we can find a simple arc $\gamma$ from $\pi(O)$ to
an arbitrary point on $\pi(\partial D^{2})$. Since the group of
homeomorphisms generated by $f$ acts freely on $D^{2}$ except at
$O$, it follows that $\pi : D^{2} \to D^{2}/f$ is a regular
branched covering (see \cite[page 49]{Maskit}). Therefore,
$\pi^{-1}(\gamma)$ is the union of $n$ disjoint simple arcs (with
the exception of their common endpoint $O$)
$\gamma_{0},\gamma_{1}^{}, \ldots ,\gamma_{n-1}$ which divide
$D^{2}$ into $n$ disjoint sectors, $A_{0},A_{1},\ldots, A_{n-1}$.
The hypothesis $\rho (f_{/\partial D^{2}}= 1/n$ implies that
$\gamma_{i} = f^{i}(\gamma_{0})$.

\begin{figure}
\begin{center}
\includegraphics[scale=0.5]{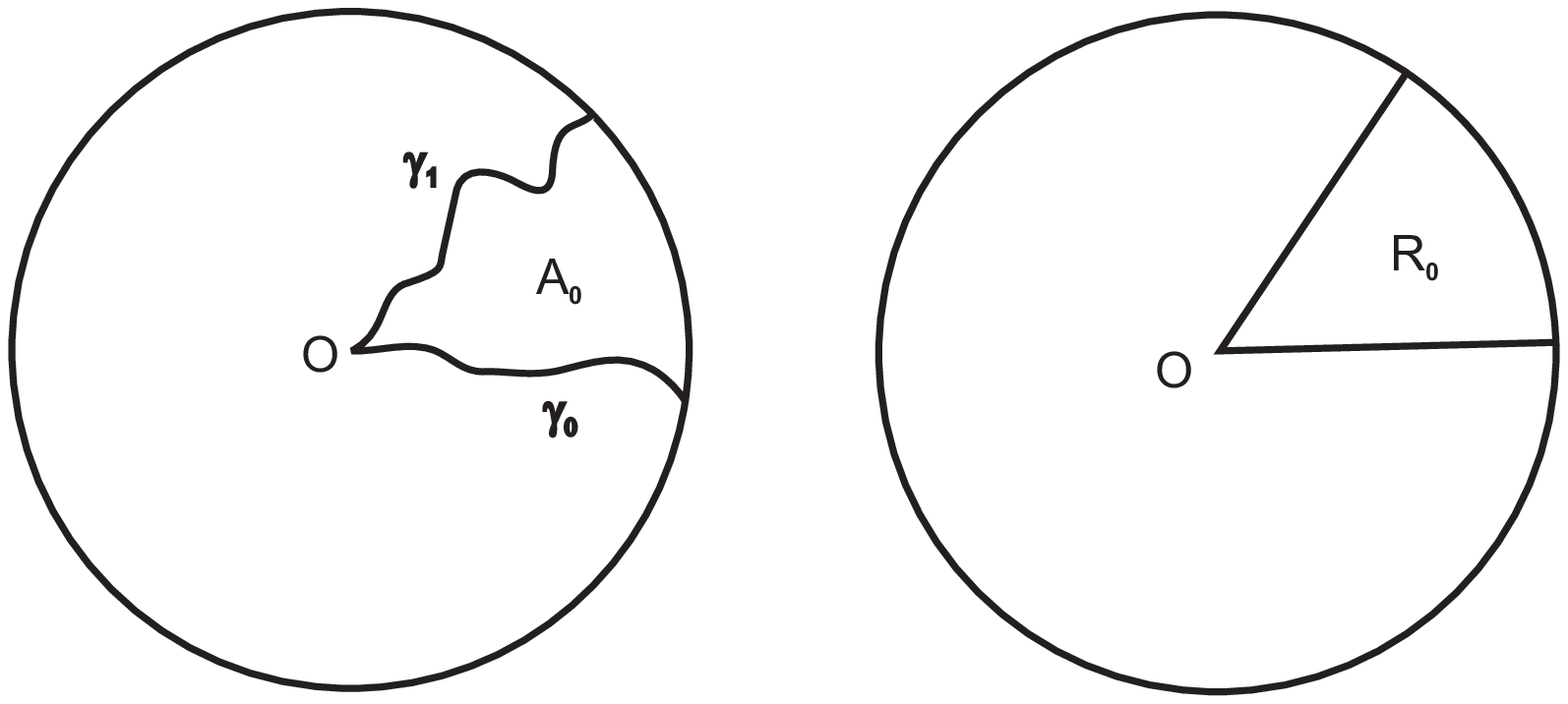}
\caption{} \label{fig3}
\end{center}
\end{figure}

Let $h$ be a homeomorphism between $A_{0}$ and $R_{0}$, the
fundamental region in $D^{2}$ of the rotation by angle $2\pi/ n$
around the origin, and such $h_{/\gamma_{1}} = rh_{/\gamma_{0}}$.
We can extend $h$ to a homeomorphism of $D^{2}$ by defining
$h_{/A_{i}}$ as $r^{i}hf^{-i}$, $r$ being the rotation of centre
$O$ and angle $2\pi / n$. It is easy t verify that $h$ is an
homeomorphism of $D^{2}$ and that $f =h^{-1}rh$.

Suppose now that $f$ is orientation-reversing. By \ref{3prop3},
$Fix(f)$ is a simple arc $\gamma$ which divides $D^{2}$ into two
topological discs $\Delta_{1}$ and $\Delta_{2}$ which are permuted
by $f$. Let $h$ be a homeomorphism between $\Delta_{1}$ and the
upper half disc $D_{1}$. We define $h$ on $\Delta_{2}$ in the
following way:
\begin{equation}
    h(y) = Sh_{/\Delta_{1}}f(y), \qquad y\in \Delta_{2},
\end{equation}

where $S$ is the reflection about the $x$-axis. It is then easy to
verify that $h$ is a homeomorphism of $D^{2}$ and this gives a
conjugacy between $f$ and $S$.
\end{proof}

\begin{rem}
Using \ref{3thm1}, it can also be shown that any periodic
homeomorphism of the annulus is topologically equivalent to an
euclidean isometry (modulo a flip of the boundary if it is not
boundary- preserving).
\end{rem}


\section{Periodic Homeomorphisms of the Sphere}

The main result of this section is:

\begin{thm}\label{4thm1}
Let $f : S^{2} \to S^{2}$ be a periodic homeomorphism. Then there
exists $r \in O(3)$ and a homeomorphism $h : S^{2} \to S^{2}$ such
that $f = hrh^{-1}$.
\end{thm}

\begin{proof} We will divide the proof of Theorem \ref{4thm1}
into two cases according to whether or not $f$ has at least one
fixed point.

Suppose first that $f$ has a fixed point. Using \ref{2lem5}, we
deduce the existence of an invariant simple closed curve $c$ which
divides $S^{2}$ into two invariant discs $D_{1}$, and $D_{2}$.

If f is orientation preserving and $f \neq Id$, then $f$ has no
fixed point on $c$ (cf. \ref{3prop2}). Therefore, by Brouwer's
fixed point theorem we know then that $f$ has at least two fixed
points; after a conjugacy, we can suppose that $f$ fixes the two
poles $N$ and $S$ of $S^{2}$. Using the results of last section,
we are able to find $n$ arcs joining $N$ and $S$ such that their
union is an invariant set under $f$. As in Section 3, we can then
construct a conjugacy between $f$ and a rotation by angle $2k\pi /
n$ around the South-North axis.

If f is orientation-reversing, then $f$ has two fixed points on
$c$. In each of the invariant disc $D_{1}$ and $D_{2}$, the fixed
point set of $f$ consists of a simple arc which joins the two
fixed points of f on c. The union of these two arcs is a simple
closed curve which coincides with the fixed point set of $f$ on
$S^{2}$. It is then easy to construct a conjugacy between $f$ and
the reflection about the equator.

Let now suppose that $f$ has no fixed point on $S^{2}$. Up to
conjugacy, we can assume that the second iterate of $f$, $f^{2}$
is a periodic rotation around the North-South axis. In particular
the points $N$ and $S$ are exchanged by $f$. For $t \in (- 1, 1)$,
let $C_{t}$ be the circle obtained by cutting the sphere by the
plane $z = t$, $D_{t}$ the disc bordered by $C_{t}$ on $S^{2}$
which contains $N$ and:
\begin{equation}
    t_{0}=\inf \{ t \in (-1,1); \quad D_{t}\cap f(D_{t}) = \emptyset \}.
\end{equation}

We write $D = D_{t_{0}}$ and $C = C_{t_{0}}$ for convenience. Then
$D$ meets $f(D)$ on its boundary and only on its boundary (see
Figure 4). Let $P_{0}\in C \cap f(C)$ and $P_{1}, P_{2},\ldots,
P_{n-1}$, the orbit of $P_{0}$ under $f$. The points $P_{0},
P_{2},\ldots, P_{n}$  and $P_{1}, P_{3},\ldots, P_{n-1}$ are
distinct because $f^{2}$ is a rotation of period $n / 2$.

Suppose that there exists $i \in \{1, 3,\ldots, n - 1 \}$ such
that and $P_{0}$ and $P_{i}=f^{i}(P_{0})$ coincide. Then $P_{0}$,
$S$ and $N$ are fixed by $f^{2i}$ so $f^{2i} = Id$. Therefore $2i
= n$. Let $b_{0}$ be the arc of great circle that joins $N$ to
$P_{0}$ in $D$ and $b_{n/2}$ its image under $f^{n/2}$. Then $b =
b_{0}\cup b_{n/2}$ is a simple arc joining $N$ and $S$ and not
meeting its first $(n/2) - 1$ iterates under $f$ away from $N$ and
$S$. These arcs divide the sphere into $n/2$ sectors and we can
build a conjugacy between $f$ and the composition of a rotation of
period $n/2$ around the North-South axis with a reflection about
the equator.

Suppose now that the points $P_{0}, P_{1},\ldots, P_{n-1}$ are
distinct. Let $b_{0}$ an arc of great circle joining $N$ and
$P_{0}$ in $D$ and $b^{\prime}_{0}$ an arc joining $S$ to $P_{0}$
in $f(D)$, disjoint from $f(b_{0})$ and from its first $n - 1$
iterates (which is possible since $f^{2}$ is a rotation). The
union of these two arcs is again a simple arc joining $N$ and $S$
which does not meet its first $n - 1$ iterates under $f$ away from
$N$ and $S$. The union of this arc and its iterates divides the
sphere $S^{2}$ into $n$ disjoint sectors. In that case, $f$ is
topologically equivalent to the composition of a rotation of
period $n$ around the North-South axis with a reflection about the
equator.
\end{proof}

\begin{figure}
\begin{center}
\includegraphics[scale=0.5]{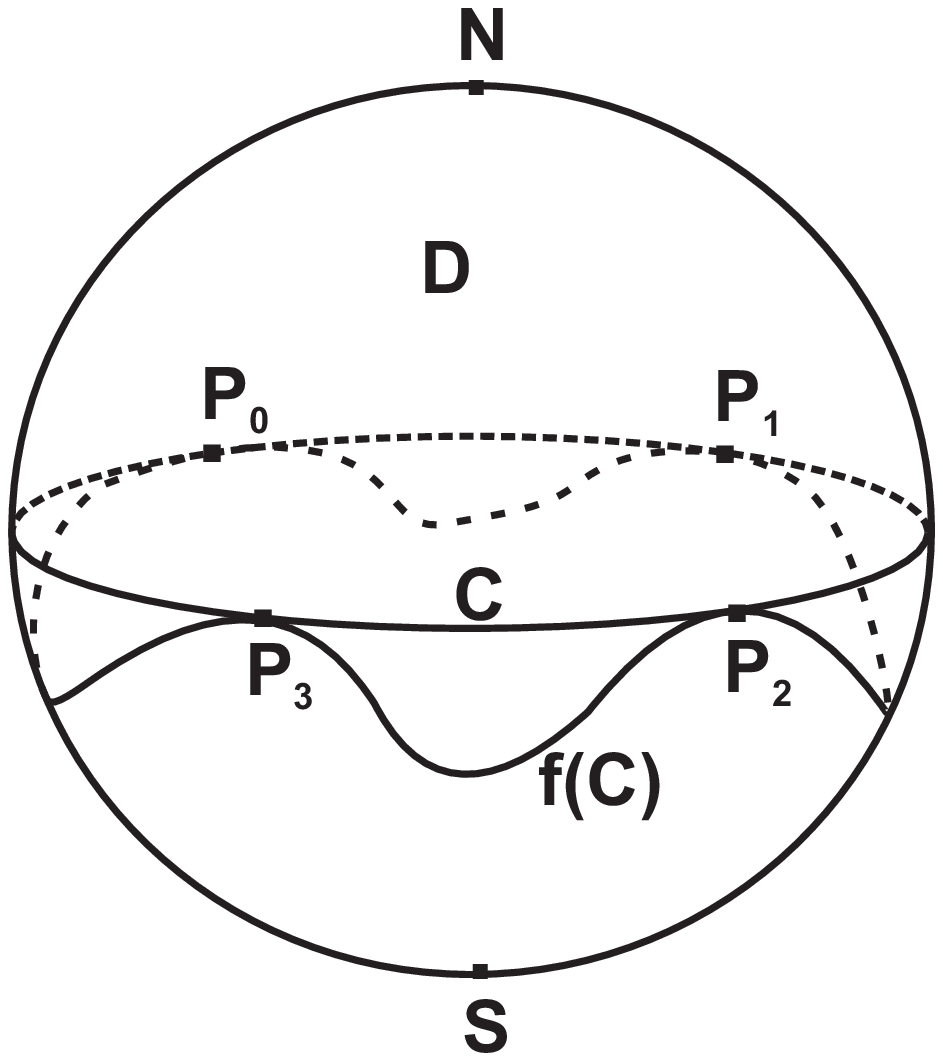}
\caption{} \label{fig4}
\end{center}
\end{figure}

\begin{cor}\label{4cor2}
Let $f : \mathbb{R}^{2} \to \mathbb{R}^{2}$ be a periodic
homeomorphism. Then $f$ is topologically conjugate to a finite
order rotation around the origin or to the reflection about the
$x$-axis.
\end{cor}

\begin{proof}
We can extend $f$ to a homeomorphism of the
Sphere $S^{2}$ by identifying the plane $\mathbb{R}^{2}$ with the
complement of the North pole using the stereographic projection.
Looking at the proof of \ref{4thm1}, $f$ is either equivalent to a
rotation around the North-South pole or to a reflection about a
great circle which we can assume to pass through the north pole
$N$. It is not difficult to show that the conjugacy can be chosen
to fix also the North pole $N$. This equivalence induces,
therefore, a topological equivalence between $f$ and a rotation or
a reflection about the $x$-axis.
\end{proof}

\begin{rem} The investigation of periodic homeomorphisms on
surfaces of positive genus has been studied extensively. We cannot
give here a complete bibliography on the subject. We would just
like to cite original works of K{\'e}r{\'e}kjart\`{o}
\cite{Kerekjarto2} and Nielsen \cite{Nielsen} which lead to the
conclusion that a periodic homeomorphism of a Riemannian surface
of positive genus is conjugate to a conformal isometry.
\end{rem}


\bibliographystyle{amsplain}
\bibliography{Periodic}
\end{document}